\documentclass[11pt,amssym,twoside]{article}
\usepackage{amssymb}
\newtheorem{thm}{Theorem}[section]
\newtheorem{lem}[thm]{Lemma}

\newtheorem{prop}[thm]{Proposition}

\newtheorem{mthm}[thm]{Main Theorem}

\newtheorem{defn}{Definition}

\def\Ab{{\mathbf A}}
\def\Ad{{\Bbb A}}

\def\C{{\Bbb C}}

\def\H{{\mathbb H}}

\def\Q{{\Bbb Q}}

\def\Z{{\Bbb Z}}
\def\Sn{{Spec(\Z[\zeta_n,1/n])}}
\def\Agn{{{\cal A}_g(n)}}
\def\Agm{{{\cal A}_g(m)}}
\def\Acgn{{{\cal A}^*_g(n)}}
\def\Acgm{{{\cal A}^*_g(m)}}
\def\Aogn{{{\cal A}^P_g(n)}}

\def\Abgn{{{\cal A}^B_g(n)}}

\def\Aocgn{{{\cal A}^{P*}_g(n)}}
\def\Aocgm{{{\cal A}^{P*}_g(m)}}

\def\Bgn{{{\cal B}_g(n)}}
\def\Bogn{{{\cal B}^P_g(n)}}

\def\a{{\alpha}}

\def\et{{\eta}}

\def\k{{\chi}}

\def\l{{\lambda}}

\def\o{{\omega}}
\def\O{{\Omega}}

\def\s{{\sigma}}

\def\z{{\zeta}}

\title{On Atkin-Lehner correspondences on Siegel spaces}

\author{Arash Rastegar}

\begin{document}

\maketitle
\begin{abstract}
We introduce a higher dimensional Atkin-Lehner theory for
Siegel-Parahoric congruence subgroups of $GSp(2g)$. Old
Siegel forms are induced by geometric correspondences on Siegel
moduli spaces which commute with almost all local Hecke algebras.
We also introduce an algorithm to get equations for moduli spaces of 
Siegel-Parahoric level structures,
once we have equations for prime levels and square prime levels 
over the level one Siegel space.
This way we give equations for an infinite tower of Siegel spaces
after N. Elkies who did the genus one case.
\end{abstract}


\subsection*{Introduction}

Classical Atkin-Lehner theory for $GL(2)$ has two conceptual
ingredients: The first one is Casselman's theory of local new-forms for $GL(2)$
[Cas], generalized to $GSp(4)$-case by R. Schmidt for square-free level
[Sch]. But, at the moment there is no general local theory of
new-forms available for $GSp(2g)$.
The second ingredient of Atkin-Lehner theory is strong multiplicity one for cuspidal
automorphic representations of $GL(2)$. 
 But, multiplicity one fails to
hold in higher $GSp(2g)$'s. Therefore, we regard all Siegel forms
with the same Hecke eigenvalues as old Siegel-forms to be old
forms again. Despite this loss, Atkin-Lehner geometric theory can be
generalized to Siegel-parahoric congruence subgroups of $GSp(2g)$.

We use elements of the Weyl group
of $GSp(2g)$ to construct new congruence subgroups sandwiched
between the Siegel-parahoric congruence group $\Gamma^P(n)$
and what we call the diagonal congruence group $\Gamma^D(n)$.
These groups are all defined in terms of the mod-$n$ reduction of
elements in $Sp(2g,\Z)$.
We shall use these congruence subgroups
to introduce geometric correspondences on Siegel moduli spaces
with explicit moduli interpretations. 
These correspondences induce an injection of a number of copies of
Siegel modular forms with Siegel-parahoric $n$-level structure inside the space of
Siegel modular forms of Siegel-parahoric level $pn$ for $(p,n)=1$. Since our
correspondences commute with all $q$-Hecke correspondences for
$(q,pn)=1$, any satisfactory local definition of $p$-old Siegel
forms will imply that our correspondences introduce a higher
dimensional Atkin-Lehner theory for $p$-old Siegel forms of
Siegel parahoric level $pn$. By a satisfactory local theory of
new-forms, we mean that any new eigenform of almost all local Hecke algebras
should be an eigenform of all local Hecke algebras prime to the
level. Also, eigenforms with eigenvalues repeated in lover levels
should be considered as old-forms. 

Next, we give algebraic equations for moduli spaces of Siegel-Parahoric level structures,
by taking fiber products of the same moduli spaces with prime power 
level structure, which in turn can be constructed via algebraic equations for
Siegel moduli spaces of Siegel-Parahoric prime level and square prime level structures.
This way we can give equations for an infinite tower of level structures
after N. Elkies who did the genus one case [Elk].

\section{Atkin-Lehner theory}

\subsection{Siegel-parahoric subgroup of $GSp(2g)$}

The Chevalley group scheme $GSp(2g)$ is defined as the set of
matrices $P\in GL(2g)$ with $^tPJP=\lambda(P)J$ where
$\lambda(P)\in GL(1)$ and $J=antidiag(I_g,-I_g)$ where $I_g$
denotes the $g\times g$ identity matrix. The representation
$\lambda$ is called the multiplier representation. Similar to
above, we also use two by two matrices whose entries are $g\times
g$ matrices to represent elements of $GSp(2g)$. The symplectic
group scheme $Sp(2g)$ is defined to be the kernel of the
multiplier representation, which is the space of transformations
on the symplectic space $\Z^{2g}$ with its standard alternating
form:
$$
<,>:\Z^{2g} \times \Z^{2g} \to \Z
\\ \hspace{.5in}<(u,v) , (z,w)>
\mapsto u.^tw-v.^tz
$$

Let $T\cong G_m^{g+1}$ denote the maximal torus in $GSp(2g)$. An
element of $T$ is a diagonal matrix
$diag(a_1,...,a_g,d_1,...,d_g)$ with $a_id_i$ equal to the
multiplier. Let $M$ denote the subgroup of $GSp(2g)$ which
respects the standard decomposition $\Z^{2g}\cong\Z^g\oplus\Z^g$.
The subgroup $M$ consists of diagonal elements in $GSp(2g)$ in
their two by two representation. These elements are of the form
$diag(A,D)$ with $A.^tD=\lambda I_g$. The subgroup of $Sp(2g,\Z)$
which fixes only the first direct summand $\Z^g\subset\Z^{2g}$ is
denoted by $U$ which is the space of all $\Z$-valued bilinear
symmetric forms. Elements of $U$ are upper-triangular with $I_g$
on diagonal entries and a symmetric matrix $B$ on the upper right
corner. The subgroup $P=M\ltimes U$ is a maximal parabolic
subgroup whose elements are zero in the lower left $g\times g$
corner. This parabolic subgroup is also called the
Siegel-parahoric subgroup. Fix the Borel subgroup contained in $P$
consisting of the matrices with $A,B,0,D$ entries with $A$
upper-triangular and $D$ lower-triangular. Weyl groups of
$GSp(2g)$ and $P$ with respect to the maximal torus $T$ are
denoted by $W_G$ and $W_P$ respectively. $W_G\cong S_g\ltimes (\pm
1)^g$ and $W_P\cong S_g$ act on diagonal matrices
$diag(a_1,...,a_g,d_1,...,d_g)$ by permutation or exchange of the
$a_i$'s and the $d_i$'s.

\subsection{Congruence subgroups of $Sp(2g,\Z)$}

A discrete subgroup $\Gamma\subset Sp(2g,\Q)$ is called a
congruence group, if it contains $\Gamma(n)=\Gamma^{id}(n)$ for some
positive integer $n$, where $\Gamma(n)$ is the kernel of
reduction map modulo $n$ on $Sp(2g,\Z)$:
$$
\Gamma(n)=\{\gamma \in Sp(2g,\Z)| \gamma \equiv I_{2g}\in
Sp(2g,\Z/n\Z) \hspace{.05in}(\textmd{mod} \hspace{.05in} n) \}.
$$
The Siegel-parahoric congruence group is defined by 
$$
\Gamma^P(n)=\{\gamma\in Sp(2g,\Z)|\gamma\equiv
\left(\begin{array}{cc} *&*\\0&* \end{array}\right), (\textmd{mod}
\hspace{.05in}n)\}
$$
and the diagonal congruence group by
$$
\Gamma^D(n)=\{\gamma\in Sp(2g,\Z)|\gamma\equiv
diag(*,...,*),(\textmd{mod} \hspace{.05in}n)\}.
$$

$\Gamma^P(n)$ and $\Gamma^D(n)$ are examples of congruence groups.
The significance of congruence groups is that they carry
arithmetic information. In this paper, we are only interested in
congruence subgroups of $Sp(2g,\Z)$.

By the theory of Tits systems [Hum], every parabolic subgroup of
$GSp(2g)$ is conjugate to a "standard" parabolic subgroup which
contains $B$. Each standard parabolic subgroup $P_I$ corresponds
to one of the $2^{g+1}$ subsets $I\subset S$ where $S$ is a
minimal generating set for the Weyl group $W_G$ consisting of
involutions $\rho_i$ which are elements of order 2. $S$ can be
taken as the set of simple reflections corresponding to the base
of the root system determined by B. In fact, $P_I=BW_I B$ where
$W_I\subseteq W_G$ is the subgroup generated by elements in $I$.
In particular, $P_{\emptyset}=B$ and $P_S=GSp(2g)$.

By Bruhat decomposition theorem, $GSp(2g)=\cup B\sigma B$ is a
decomposition to disjoint subsets where $\sigma$ runs in the Weyl
group $W_G$. Two such double cosets coincide if and only if the
middle Weyl elements coincide. One can assume that axioms of Tits
systems are satisfied. Namely,

(T1) For $\rho\in S$ and $\sigma\in W_G$ we have $\rho
B\sigma\subset B\sigma B \cup B\rho\sigma B$.

(T2) For $\rho\in S$ we have $\rho B\rho \neq B$.
\\An expression $\sigma=\rho_1...\rho_k$ with $\rho_i\in S$ is
called reduced if $k$ is as small as possible. the minimal length
$k$ of a reduced expression is denoted by $\ell(\sigma)$. By
convention, $\ell(\sigma)=0$ if and only if $\sigma=id$ and
$\ell(\sigma)=1$ if and only if $\sigma\in S$. Tits axioms imply
that for $\rho\in S$ we have $\ell(\rho\sigma)=\ell(\sigma)\pm 1$.
This implies that for the reduced form $\sigma=\rho_1...\rho_k$
and $I=\{\rho_1,...,\rho_k\}$, the parabolic subgroup $P_I$ is
generated $B$ and $\sigma B\sigma^{-1}$ or by $B$ and $\sigma$.
Given a choice of a Borel subgroup $S$ is precisely the set of
elements in $W_G$ such that $B\cup B\sigma B$ is a group. $P_I$ is
conjugate to $P_J$ implies that $P_I=P_J$ and $P_I\subset P_J$
implies that $I\subset J$.

We define a parabolic congruence subgroup $\Gamma^{P_I}(n)\subset
Sp(2g,\Z)$ to be the set of elements which reduce to $P_I$ modulo
$n$:
$$
\Gamma^{P_I}(n)=\{\gamma \in Sp(2g,\Z)| \gamma \in P_I\subset
GSp(2g,\Z/n\Z) \hspace{.05in}(\textmd{mod} \hspace{.05in} n) \}.
$$

Fix a generating set $S$ for $W_G$ consisting of $g-1$ pairs in
$S_g\subset W_g$ and a nonzero element in $(\pm 1)^g$. Now, there
exists a generating set $w_1,...,w_g$ for the weyl group such that
$w_i$ have increasing length if represented in reduced form in
term of involutions in $S$.

Assume that $w_{k+1}=\rho_k w_k$for all $k$ where $\rho_k$ is an
involution. Also assume $w_1$,...,$w_{g-1}$ generate $S_g$. We
have $\Gamma^{P_{I_{g-1}}}(n)=\Gamma^P(n)$ where $P$ is the
maximal parabolic we fixed in notations. The fact that $B$
together with $I_k=\{w_1,...,w_k\}$ can not generate any of
$w_{k+1}$,...,$w_g$ implies that each $P_{I_k}$ has exactly
$\frac{g!}{(k+1)!}2^g$ conjugates of the form $\sigma
P_{I_k}\sigma^{-1}$ for $\sigma$ in $W_G$. To these we associate
$\frac{g!}{(k+1)!}2^g$ parabolic congruence subgroups of the form
$\sigma \Gamma^{P_{I_k}}(n)\sigma^{-1}$ for $\sigma$ in $W_G$. In
particular, $\Gamma^B(n)$ has $g! 2^g$ conjugates in $Sp(2g,\Z)$
and each $\Gamma^{P_{I_k}}(n)$ has $k+2$ conjugates in
$\Gamma^{P_{I_{k+1}}}(n)$.

We have made a nested family of parabolic congruence groups
contained in the maximal parabolic congruence group $\Gamma^P(n)$
consisting of levels $0$ to $g-1$. The $k$-th level is formed by
$\frac{g!}{(k+1)!}2^g$ congruence groups. Each group in level $k$
contains $k+1$ groups in level $k-1$ for $k=1$ to $g-1$. There
are $2^g$ congruence groups in level $g-1$ and the congruence
groups in level $0$ are the $g!2^g$ conjugates by elements in
$W_g$ of $\Gamma^B(p)$ which are lying inside $Sp(2g,\Z)$.

\subsection{Siegel moduli spaces}

A general reference for the arithmetic of Siegel moduli spaces is
[Fal-Cha]. Let $A_g$ denote the moduli stack of principally
polarized abelian schemes of relative dimension $g$. By a
symplectic principal level-$n$ structure, we mean a symplectic
isomorphism $\a:A[n]\to(\Z/n\Z)^{2g}$, where $(\Z/n\Z)^{2g}$ is
equipped with the standard non-degenerate skew-symmetric pairing.

Let $\z_n$ denote an $n$-th root of unity for $n \geq 3$. The
moduli scheme classifying the principally polarized abelian
schemes over $\Sn$ together with a symplectic principal level-$n$
structure is a scheme over $\Sn$ and will be denoted by $\Agn$.
The symplectic group $Sp(2g,\Z/n\Z)$ acts on $\Agn$ as a group of
symmetries by acting on level structures. We will recognize these
moduli spaces and their etale quotients under the action of
subgroups of $Sp(2g,\Z/n\Z)$ as Siegel spaces.

A $\Gamma^P(n)$-level structure of type I on $(A,\l)$ is choice
of a subgroup $H\subset A[n]$ of order $n^g$ which is totally
isotropic with respect to the Weil pairing induced by $\l$. A
$\Gamma^P(n)$-level structure of type II on $(A,\l)$ is choice of
a principally polarized isogeny $(A_1,\l_1)\to (A_2,\l_2)$ of
degree $n^g$. By a principally polarized isogeny, we mean an
isogeny $\s :A_1 \to A_2$ such that $\s\circ \l_2\circ \s^t\circ
\l_1^{-1}$ is multiplication by an integer. For $n\geq 3$ type I
and type II $\Gamma^P(n)$-level structures induce isomorphic
moduli schemes over $Spec(\Z [1/n])$ [DJo]. We denote this moduli
scheme by $\Aogn$. There exists a natural involution
$$
w^P_n:\Aogn \to \Aogn
$$
taking $(\s :(A_1,\l_1)\to (A_2,\l_2))$ to $(\s^t
:(A_2,(\l_2^t)^{-1}) \to (A_1,(\l_1^t)^{-1}))$ which we call the
Atkin-Lehner involution.

A $\Gamma^B(n)$-level structure of type I on $(A,\l)$ is choice
of $g$ subgroups $H_i\subset A[n]$ of order $n^i$ with $H_1
\subset ... \subset H_g$ where $H_g$ is totally isotropic. A
$\Gamma^B(n)$-level structure of type II on $(A,\l)$ is choice of
a chain of $g$ isogenies $(A_0,\l_0) \buildrel \a \over \to ...
\buildrel \a \over\to (A_g,\l_g)$ each of degree $n$ which
satisfy $n.id_{A_i}=\a^i\circ \l_0^{-1}\circ (\a^t)^g\circ \l_g
\circ \a^{g-i}$ for all $i=1,...,g$. In case $n\geq 3$ type I and
type II $\Gamma^B(n)$-level structures induce isomorphic moduli
schemes over $Spec(\Z [1/n])$ [DJo]. We denote this moduli scheme
by ${\cal{A}}^B_{g}(n)$. There also exists a natural involution
$$
w^B_n:\Abgn\to \Abgn
$$
taking $ ((A_0,\l_0) \buildrel \a \over \to ... \buildrel \a \over
\to (A_g,\l_g))$ to $((A_g,(\l_g^t)^{-1}) \buildrel \a^t \over \to
... \buildrel \a^t \over \to (A_0,(\l_0^t)^{-1}))$ which commutes
with the Atkin-Lehner involution under the natural projection
between the Siegel spaces.
$$
\begin{array}{ccc}
\Abgn & \buildrel {w^B_n} \over \longrightarrow & \Abgn\\
\downarrow &                     & \downarrow \\
\Aogn & \buildrel {w^P_n} \over  \longrightarrow & \Aogn
\end{array}
$$

A $\Gamma^T(n)$-level structure on $(A,\l)$ is choice of $2g$
subgroups $H_i\subset A[n]$ for $i=1$ to $2g$, each isomorphic to
$(\Z/n\Z)$ such that $H_1\oplus ...\oplus H_g$ and $H_{g+1}\oplus
...\oplus H_{2g}$ are totally isotropic subgroups of order $n^g$
which do not intersect with $H_i\oplus H_{g+i}$ hyperbolic for
$i=1$ to $g$. For $A$ and $A'$ abelian schemes over the schemes
$S$ and $S'$ respectively, we define a morphism from $(S,A,\l
,H_1,...,H_{2g})$ to $(S',A',\l' ,H'_1,...,H'_{2g})$ to be a pair
of morphisms $(f,g)$ where $f:S\to S'$ and $g:A\to A'$ satisfy
$g^*(\l')=\l$ and $g(H_i)=H'_i$ for all $1\leq i \leq 2g$. Also
we want the pair $(f,g)$ to induce an isomorphism $A \simeq
S\times_{S'}A'$. Having these morphisms defined, we have formed a
category $\Ab ^T_g(n)$. The functor $\pi :\Ab ^T_g(n)\to Sch$
defined by $\pi (S,A,\l ,H_1,...,H_{2g})=S$ makes $\mathbf A ^T
_g(n)$ into a stack in groupoids over $S$. The $1$-morphism of
stacks $\pi':\mathbf A ^T_g(n)\to \mathbf A_g$ defined by
$\pi'(S,A,\l ,H_1,...,H_{2g})=(S,A,\l )$ is representable and is
a proper surjective morphism. For $n\geq 3$ we get a separated
scheme of finite type $A^T_g(n)$ which is smooth over
$Spec(\Z[1/n])$.

Let $\H_g$ denote the Siegel upper half-space, which consists of
the set of complex symmetric $g \times g$ matrices $\O$ with
$\Im(\O)$ positive definite. As a complex analytic stack ${\cal
{A}}_{g/\C}$ is the quotient of Siegel upper half-space $\H_g$ by
the action of $Sp(2g,\Z)$ via M\"obius transformations. The
family of principally polarized abelian varieties over $\H_g$ is
given by $A(\O)=\C^g/(\Z^g \oplus \O.\Z^g)$. To any congruence
subgroup $\Gamma\subset Sp(2g,\Z)$ one can associate the quotient
$\Gamma\backslash \H_g$ which is a Siegel moduli-space with some
extra level structure. The corresponding level structure can be
made explicit. Indeed, $\Gamma(n)\backslash\H_g$ corresponds
to $\Agn_{/\C}$ whose quotient under the action of the symplectic
group $Sp(2g,\Z/n\Z)$ is ${\cal {A}}_{g/\C}$. Any congruence
subgroup $\Gamma(n)\subset \Gamma \subset Sp(2g,\Z)$
corresponds to the quotient of $\Agn_{/\C}$ under the action of
the finite group $\Gamma/\Gamma(n)$. This helps to associate
explicit level structures to the space $\Gamma\backslash\H_g$
which makes it a moduli space.

\subsection{The $\Gamma^{P_I}(p)$-level structure}

The Atkin-Lehner involution can be easily generalized from
$GL(2)$ [Atk-Leh] to the higher dimensional case $GSp(2g)$. This
generalization involves Siegel-parahoric and Borel congruence
groups. On the other hand, local considerations show that $p$-old
Siegel modular forms with respect to the Siegel-parahoric or Borel
congruence groups of level $pn$ contain several copies of Siegel
forms of level $n$. This implies that a single Atkin-Lehner
involution would not do the job of geometrically generating the
$p$-old part. In this section,  we intend to introduce geometric
correspondences which complement the role of Atkin-Lehner
involution.

Let $p$ be a prime not dividing the positive integer $n$ and let
${\cal {A}}^T_g(p)_{/\C}$ and ${\cal {A}}^{T,P}_g(p,n)_{/\C}$
denote the Siegel spaces associated to congruence groups
$\Gamma^T(p)$ and $\Gamma^{T,P}(p,n)=\Gamma^T(p)\cap
\Gamma^{P}(n)$, respectively. The group $\Gamma^T(p)$ remains
invariant under conjugation by elements in $(\pm 1)^g\subset
W_g$. So $(\pm 1)^g$ acts on ${\cal {A}}^T_g(p)_{/\C}$ and ${\cal
{A}}^{T,P}_g(p,n)_{/\C}$ by $2^g$ involutions. Let ${\cal
{A}}^{P_{I_k}}_g(p)_{/\C}$ and ${\cal
{A}}^{P_{I_k},P}_g(p,n)_{/\C}$ denote the Siegel spaces
associated to the congruence groups $\Gamma^{P_{I_k}}(p)$ and
$\Gamma^{P_{I_k},P}(p,n)=\Gamma^{P_{I_k}}(p)\cap \Gamma^{P}(n)$
respectively. We have a chain of etale maps
$$
{\cal {A}}_g^{B,P}(p,n) \rightarrow ... \rightarrow {\cal
{A}}^{P_{I_k},P}_g(p,n) \rightarrow {\cal
{A}}^{P_{I_{k+1}},P}_g(p,n) \rightarrow ... \rightarrow {\cal
{A}}_g^{P}(pn).
$$
Since each congruence group $\Gamma^{P_{I_{k+1}}}(p)$ on the
$(k+1)$-th level contains $k+2$ conjugates (by Weyl elements) of
$\Gamma^{P_{I_k}}(p)$
 on the $k$-th level, we expect that for each $k$
we get $k+2$ copies of forms on ${\cal
{A}}^{P_{I_{k+1}},P}_g(p,n)$ injecting in forms on ${\cal
{A}}^{P_{I_k},P}_g(p,n)$. We will use the geometry of ${\cal
{A}}_g^{B,P}(p,n)$ to give a geometric construction of these $k+2$
copies.

In order to simplify the notations, let us forget the
$\Gamma^{P}(n)$-level structure which is auxiliary. We get a
chain of etale maps
$$
{\cal {A}}_g^{B}(p) \rightarrow ... \rightarrow {\cal
{A}}^{P_{I_k}}_g(p) \rightarrow {\cal {A}}^{P_{I_{k+1}}}_g(p)
\rightarrow ... \rightarrow {\cal {A}}_g^{P}(p)
$$
which corresponds to a chain of congruence groups
$$
\Gamma^{B}(p) \hookrightarrow ... \hookrightarrow
\Gamma^{P_{I_k}}(p) \hookrightarrow \Gamma^{P_{I_{k+1}}}(p)
\hookrightarrow ... \hookrightarrow \Gamma^{P}(p).
$$
Each $\Gamma^{P_{I_k}}(p)$ maps to $\Gamma^{P_{I_{k+1}}}(p)$ by
$k+1$ maps: natural inclusion and conjugation by representatives
$\sigma_{k+1}\in W_k\backslash W_{k+1}$ followed by inclusion,
where $W_k$ is the subgroup of the Weyl group generated by
$w_1,...,w_k$. Inclusion induces the natural projection map
$\pi_{P_{I_k},P_{I_{k+1}}}:{\cal {A}}^{P_{I_k}}_g(p) \rightarrow
{\cal {A}}^{P_{I_{k+1}}}_g(p)$. Conjugation by $\sigma_{k+1}$
induces an inclusion
$$
\sigma_{k+1}\Gamma^{P_{I_k}}(p)\sigma_{k+1}^{-1}=
\Gamma^{\sigma_{k+1}P_{I_k}\sigma_{k+1}^{-1}}(p)\hookrightarrow
\Gamma^{P_{I_{k+1}}}(p).
$$
This inclusion corresponds to another projection from a different
moduli-space
$$
{\cal {A}}^{\sigma_{k+1}P_{I_k}\sigma_{k+1}^{-1}}_g(p)
\rightarrow {\cal {A}}^{P_{I_{k+1}}}_g(p).
$$
Conjugation by $\sigma_{k+1}$ identifies ${\cal
{A}}^{P_{I_k}}_g(p)$ with ${\cal
{A}}^{\sigma_{k+1}P_{I_k}\sigma_{k+1}^{-1}}_g(p)$. The
moduli-space ${\cal {A}}^T_g(p)$ is the appropriate moduli space
to geometrically realize all the endomorphisms
$$
v^{\sigma}_p:{\cal {A}}^T_g(p) \to {\cal {A}}^T_g(p)
$$
induced by conjugation via elements $\sigma$ in the Weyl group
$W_G$. In fact, the following diagrams are commutative

$$
\begin{array}{ccc}
{\cal {A}}_g^T(p) & \buildrel {v^{\sigma}_p} \over \longrightarrow
& {\cal {A}}_g^T(p)\\
\downarrow &                     & \downarrow \\
{\cal {A}}_g^{P_I}(p) & \buildrel {\sigma .} \over \longrightarrow
& {\cal {A}}_g^{\sigma P_I \sigma^{-1}}(p) \\
\downarrow &                     & \downarrow \\
{\cal {A}}_g^{P_J}(p) & \buildrel {\sigma .} \over \longrightarrow
& {\cal {A}}_g^{\sigma P_J \sigma^{-1}}(p)
\end{array}
$$

\subsection{The geometry of Siegel spaces ${\cal {A}}_g^{B}(p)$ and ${\cal {A}}_g^{P_I,P}(p)$}

In this section we try to geometrically characterize fibers of natural
maps between moduli spaces, we have already introduced.

\begin{defn} We say that two points $x$ and $y$ on ${\cal
{A}}_g^{P_I,P}(p,n)$ are $\sigma$-connected, for an element
$\sigma$ in the Weyl group, if there exists a chain of points
$x=x_1,...,x_t=y$ on ${\cal {A}}_g^{P_I,P}(p,n)$ such that for
each $i$ there are points $x'_i$ and $x'_{i+1}$ on ${\cal
{A}}^{T,P}_g(p,n)$ mapping to $x_i$ and $x_{i+1}$ respectively,
with $x_{i+1}=v^{\sigma}_p(x_i)$ where
$$
v^{\sigma}_p:{\cal {A}}^{T,P}_g(p,n) \to {\cal {A}}^{T,P}_g(p,n)
$$
is the endomorphism induced by the action of $\sigma$ on
$p$-level structure.
\end{defn}

\begin{prop} Every fiber of the map $\pi_k:{\cal {A}}^{P_{I_k},P}_g(p,n)
\rightarrow {\cal{A}}^{P_{I_{k+1}},P}_g(p,n)$ is an equivalence
class of $\rho_{k+1}$-connected points.
\end{prop}
$\textbf {Proof}$. The $n$-level structure is auxiliary. Let $x$
be a point on the Siegel upper half-plane $\H_g$ and let $[x]^T$
denote the equivalence class containing $x$ defined by left
quotient of the Siegel upper space by $\Gamma^T(p)$. We have
$\sigma .[y]^T=[\sigma.y]^T$. The group $\Gamma^{P_{I_{k+1}}}(p)$
is generated $\Gamma^{P_{I_k}}(p)$ and $\rho_{k+1}$. So every
element in $\Gamma^{P_{I_{k+1}}}(p)$ can be written as a product
of elements of the form $\gamma_i \rho_{k+1}$ with $\gamma_i \in
\Gamma^{P_{I_k}}(p)$. Define the equivalence class
$[x]^{P_{I_k}}$ similarly. The classes $[x]^{P_{I_k}}$ and
$[\gamma\rho_{k+1}.x]^{P_{I_k}}$ are $\rho_{k+1}$-connected. So
the equivalence class $[x]^{P_{I_{k+1}}}$ is obtained by joining
the $\rho_{k+1}$-connected points.$\square$

\begin{defn} For a subset $W\subset W_g$, we say that
two points $x$ and $y$ on ${\cal{A}}_g^{P_I,P}(p,n)$ are
$W$-connected, if there exists a chain of points
$x=x_1,...,x_t=y$ on ${\cal{A}}_g^{P_I,P}(p,n)$, such that for
each $i$, $x_i$ and $x_{i+1}$ are $\sigma$-connected for some
$\sigma\in W$.
\end{defn}

\begin{prop} Every fiber of the map ${\cal{A}}_g^{B,P}(p,n)\to
{\cal{A}}_g^{P_{I_k},P}(p,n)$ is an equivalence class of
$W_k$-connected points.
\end{prop}

\begin{prop} Every fiber of the map ${\cal{A}}_g^P(pn)\to
{\cal{A}}_g^P(n)$ is an equivalence class of $\sigma$-connected
points for some nonzero representative $\sigma$ of $W_P \backslash
W_G$.
\end{prop}
$\textbf {Proof}$. The Bruhat-Tits decomposition modulo $p$
implies that the congruence groups $\Gamma^P(p)$ and
$\sigma\Gamma^P(p)\sigma^{-1}$ generate $Sp(2g,\Z)$. This is a
consequence of the fact that any two conjugates of a maximal
parabolic subgroup over $\Bbb {F}_p$ generate the whole algebraic
group $GSp(2g,\Bbb {F}_p)$. $\square$

\subsection{$p$-old Siegel modular forms on ${\cal {A}}_g^{P}(p)$}

Let $\Bogn$ denote the universal abelian variety over $\Aogn$.
The Hodge bundle $\o$ is defined to be the pull back via the zero
section $i_0:\Aogn \to \Bogn$ of the line bundle
$\wedge^{top}\Omega^1_{\Bgn/\Agn}$. The Hodge bundle is an ample
invertible sheaf on $\Aogn$. Let $R$ be a $\Z [1/n]$-module. By a
$\Gamma^P(n)$-Siegel modular form of weight $k$ with coefficients
in $R$, we mean an element in $H^0(\Aogn /R,\o^{\otimes k}/R)$.
The same notation is used for other congruence subgroups, but in
this paper we focus on Siegel modular forms with respect to
Siegel parahoric congruence subgroup $P$ and Borel congruence
subgroup $B$.

If we pull back the Hodge bundle $\o$ to the Siegel upper
half-space, the pull back is canonically isomorphic to
${\cal{O}}_{\H_g} \otimes_{\C} \wedge^g(\C^g)$. A complex modular
form of weight $k$ becomes an expression of the form
$f(\O).(dz_1\wedge...\wedge dz_g)^{\otimes k}$ where $f$ is an
$\Gamma^P(n)$-invariant complex holomorphic function on $\H_g$
which is holomorphic at $\infty$. For genus $\geq 2$ the
condition, holomorphic at infinity, is automatically satisfied by
Koecher principle. Trivializing $\o$ on $\H_g$, complex modular
forms of weight $k$ are identified with holomorphic functions
$f(\O)$ on $\H_g$, satisfying the transformation rule
$f|[\gamma]_k=f$ for all $\gamma \in \Gamma^P(n)$ where
$$
f|[\gamma]_k(\O)=\et(\gamma)^{gk/2}det(C.\O+D)^{-k}f(\gamma.\O).
$$

Let $l$ be a prime not dividing $pn$. To any Siegel modular form
$f$ of weight $k$ and level $n$, one associates an irreducible
admissible representation $\pi=\bigotimes \pi_{\ell}$ of
$GSp(2g,\Ad_f)$ over $\Q_l$ [Asg-Sch]. This association is not unique,
but we use it as a motivation to understand the notion of $p$-old
form. Let $U$ denote the open subgroup of $GSp(2g,\Ad_f)$
associated to the congruence group $\Gamma^P(n)$. If $(p,n)=1$ and
$\pi^U\neq 0$, then $\pi_p$ is spherical, and it is the unique
irreducible subquotient of some unramified principal series
representation $\pi_{\k}$ with respect to Borel subgroup
$B(\Q_p)$. One can show that $\pi_{\k}^{GSp(2g,\Z_p)}$ is
one-dimensional. So the number of copies of modular forms with
respect to $\Gamma^P(n)$ inside modular forms with respect to
$\Gamma^P(pn)$ is equal to the dimension of
$\pi_{\k}^{\Gamma^P(p)}$. The mod-$p$ Bruhat-Tits decomposition
implies that, we have the following decomposition
$$
GSp(2g,\Q_p)=\coprod_{\sigma\in (\Z/2\Z)^g \subset W_g}
B(\Q_p)\sigma\Gamma^P(p).
$$
So to specify $f\in \pi_{\k}^{\Gamma^P(p)}$ it is enough to
specify it on elements $\sigma\in (\Z/2\Z)^g$. Because
$W_{P}\simeq S_g$ and the subgroup $(\Z/2\Z)^g \subset W_G$ is a
complete set of representatives for $W_{P}\backslash W_G$.
Therefore, the space of $p$-old forms on ${\cal {A}}_g^P(pn)$
consists of $2^g$ copies of forms on ${\cal {A}}_g^P(n)$. The
vector space $\pi_{\k}^{\Gamma^P(p)}$ has a basis consisting of
functions $f_1,...,f_{2^g}$ where $f_i$ is supported on
$B(\Q_p)\sigma_i\Gamma^P(p)$ and $f_i(\sigma_i)=1$. The group of
involutions $(\Z/2\Z)^g\subset W_G$ acts on the space of $p$-old
forms by $f(z)\mapsto f^{\sigma}(z):=f(\sigma .z)$ for $\sigma\in
(\Z/2\Z)^g$.

Similar considerations show that we expect $g!2^g$ copies of forms
on ${\cal {A}}_g^P(n)$ inside the space of $p$-old forms on ${\cal
{A}}_g^{B,P}(p,n)$. Following Atkin-Lehner theory, we need a
geometric characterization of the space of $p$-old forms.

\subsection{Atkin-Lehner correspondences}

Let $\pi^{}_{T,P_I}:{\cal {A}}_g^T(p)\to {\cal {A}}_g^{P_I}(p)$
and $\pi^{}_{P_I,P_J}:{\cal {A}}_g^{P_I}(p)\to {\cal
{A}}_g^{P_J}(p)$ denote the natural projection maps induced by
inclusions of the corresponding congruence groups. Pulling back
forms from level $n$ to level $np$ using the natural projection
map
$$
\pi_n:{\cal {A}}_g^{B,P}(p,n)\rightarrow {\cal {A}}_g^P(n)
$$
induces the first copy of $p$-old forms in $H^0({\cal
{A}}_g^{B,P}(p,n) ,\o^{\otimes k})$. For simplicity, let us forget
the auxiliary level structure and consider the projection
$$
\pi_1:{\cal {A}}_g^{B}(p)\rightarrow {\cal {A}}_g
$$
and $p$-old forms in $H^0({\cal {A}}_g^B(p) ,\o^{\otimes k})$. At
first glance, it seems that geometric correspondences of the form
$D^{\sigma}_B(p)=\pi^{}_{T,B*} v^{\sigma *}_p\pi^{*}_{T,B}$ should
induce more copies of $p$-old Siegel forms on ${\cal {A}}_g^B(p)$
out of the pull-back copy.
$$
\begin{array}{ccc}
{\cal {A}}_g^T(p) & \buildrel {v^{\sigma}_p} \over \longrightarrow
& {\cal {A}}_g^T(p)\\
\downarrow &                     & \downarrow \\
{\cal {A}}_g^B(p) &                         & {\cal {A}}_g^B(p)
\end{array}
$$
But $\pi_{T,B}\circ \pi_1$ commutes with $v^{\sigma}_p$ for all
$\sigma\in W_G$. Therefore, correspondences of above type
generate the same copy as the pull back copy. To disturb the
symmetry of the picture, we use Atkin-Lehner involution.
$\pi_{T,B}\circ w^B_p \circ \pi_1$ no longer commutes with
$v^{\sigma}_p$ and we can hope that correspondences of the form
$C^{\sigma}_P(p)=\pi^{}_{T,B*} v^{\sigma
*}_p\pi^{*}_{T,B}w^{B*}_p$ could generate more copies of $p$-old
forms.
$$
\begin{array}{ccccccc}
& & {\cal {A}}_g^T(p) & \buildrel {v^{\sigma}_p} \over
\longrightarrow
& {\cal {A}}_g^T(p)& & \\
& & \downarrow &                     & \downarrow & & \\
& &{\cal {A}}_g^B(p) & & {\cal {A}}_g^B(p) & \buildrel {w^B_p}
\over \longrightarrow & {\cal {A}}_g^B(p)
\end{array}
$$
To generate $p$-old forms in $H^0({\cal {A}}_g^B(p) ,\o^{\otimes
k})$ we should use $g!2^g$ correspondences $C^{\sigma}_B(p)$ for
$\sigma\in W_G$.

\begin{mthm}
The linear subspaces of $H^0({\cal {A}}_g^{B,P}(p,n) ,\o^{\otimes
k})$ generated by Atkin-Lehner correspondences
$C^{\sigma}_B(p)\pi^*_n$ where $C^{\sigma}_B(p)$ is defined by
$\pi^{}_{T,B*} v^{\sigma *}_p\pi^{*}_{T,B}w^{P*}_p$
$$
\begin{array}{ccccccc}
& & {\cal {A}}_g^{T,P}(p,n) & \buildrel {v^{\sigma}_p} \over
\longrightarrow
& {\cal {A}}_g^{T,P}(p,n)& & \\
& & \downarrow &                     & \downarrow & & \\
& &{\cal {A}}_g^{B,P}(p,n) & & {\cal {A}}_g^{B,P}(p,n) & \buildrel
{w^B_p} \over \longrightarrow & {\cal {A}}_g^{B,P}(p,n)
\end{array}
$$
for $\sigma$ varying in $W_G$ give $g!2^g$ linearly independent
copies of $H^0({\cal {A}}_g^P(n) ,\o^{\otimes k})$ inside $p$-old
forms of level $pn$ living on ${\cal {A}}_g^{B,P}(p,n)$.
\end{mthm}

The corresponding theorem for $H^0({\cal {A}}_g^P(p) ,\o^{\otimes
k})$ can also be proved. However, in order to generate $p$-old
forms in $H^0({\cal {A}}_g^P(p) ,\o^{\otimes k})$ we should
divide the space of $p$-old forms on ${\cal {A}}_g^B(p)$ by the
action of $(\Z/2\Z)^{2g}\subset W_G$.

Main theorem is proved in a few stages. In the first stage, we
prove that Atkin-Lehner involution induces a second copy of
level-$n$ forms inside $p$-old part of level-$np$ forms which has
trivial intersection with the pull-back copy.
\begin{prop}
$\pi_{B,P_I}^*H^0({\cal {A}}_g^{P_I,P}(p,n) ,\o^{\otimes k})$ and
$w^{B*}_p\pi_{B,P_I}^*H^0({\cal {A}}_g^{P_I,P}(p,n) ,\o^{\otimes
k})$ as subspaces of $H^0({\cal {A}}_g^{B,P}(p,n) ,\o^{\otimes
k})$ have trivial intersection.
\end{prop}
\textbf{Proof.} Let $f$ and $g$ be nonzero Siegel modular forms in
$\pi_{B,P_I}^*H^0({\cal {A}}_g^{P_I,P}(p,n) ,\o^{\otimes k})$ with
$w^{B*}_pf=g$. Since $w^{B*}_p$ is an involution, $f\pm g$ are
eigenforms of $w^{B*}_p$ and the proposition follows from the
following
\begin{lem}
Any Siegel form which is eigenform of $w^{B*}_p$ on ${\cal
{A}}_g^{B,P}(p,n)$ and also pull back of a Siegel form on ${\cal
{A}}_g^{P_I,P}(p,n)$ vanishes if $P_I\subsetneqq B$.
\end{lem}
\textbf{Proof.} The zero locus of such an eigenform is pull back
of the zero locus of a form living on ${\cal {A}}_g^{P_I,P}(p,n)$
and also $w^{B*}_p$-invariant. This contradicts density of Hecke
orbit [Cha]. $\square$

In the second stage, we show that Atkin-Lehner correspondences
induce $g!2^g$ non-intersecting copies of level-$n$ forms inside
$p$-old part of level-$np$ forms.

\begin{lem} Let $D^{\sigma}_B(p)=\pi^{}_{T,B_*} v^{\sigma
*}_p\pi^{*}_{T,B}$ where $\pi^{}_{T,B}:{\cal {A}}_g^{T,P}(p,n) \to
{\cal {A}}_g^{B,P}(p,n)$ is the natural projection. Then, for
$\sigma,\sigma'\in W_G$ the correspondences
$D^{\sigma}_B(p)D^{\sigma'}_B(p)$ and $D^{\sigma\sigma'}_B(p)$
acting on any linear subspace of $H^0({\cal {A}}_g^{B,P}(p,n)
,\o^{\otimes k})$ generate the same image subspaces.
\end{lem}
\textbf{Proof.} For simplicity, let us forget the auxiliary
$n$-level structure. Then, lemma follows from commutativity of
the following diagram
$$
\begin{array}{ccccccc}
{\cal {A}}_g^{T}(p)& \buildrel {v^{\sigma}_p} \over
\longrightarrow & {\cal {A}}_g^{T}(p) & \buildrel {v^{\sigma'}_p}
\over \longrightarrow
& {\cal {A}}_g^{T}(p)& & \\
\downarrow & & \downarrow &  & \downarrow & & \\
{\cal {A}}_g^{\sigma\sigma' B\sigma'^{-1}\sigma^{-1}}(p) &
\buildrel {\sigma} \over \longrightarrow & {\cal {A}}_g^{\sigma'
B\sigma'^{-1}}(p) & \buildrel {\sigma'} \over \longrightarrow &
{\cal {A}}_g^B(p) & &
\end{array}
$$
and, $v^{\sigma}_p\circ v^{\sigma'}_p=v^{\sigma\sigma'}_p $ and
$\sigma\circ \sigma'=\sigma\sigma'$ hold for $\sigma,\sigma'\in
W_G$.$\square$

\begin{prop}
Every two linear subspaces of $H^0({\cal {A}}_g^{B,P}(p,n)
,\o^{\otimes k})$ generated by correspondences
$C^{\sigma}_B(p)\pi^*_n$ acting on $H^0({\cal {A}}_g^{P}(n)
,\o^{\otimes k})$ for $\sigma\in W_G$ have trivial intersection.
\end{prop}
\textbf {Proof.} By the previous lemma, it is enough to show that
the linear subspaces generated by correspondences $\pi^{}_{T,B_*}
v^{\sigma *}_p\pi^{*}_{T,B}w^{B*}_p\pi_1^*$ and $w^{B*}_p\pi_1^*$
have such intersection. Suppose $\pi^{}_{T,B_*} v^{\sigma
*}_p\pi^{*}_{T,B}f=g$ for nonzero Siegel modular forms on ${\cal
{A}}_g^B(p)$ which are elements of $w^{B*}_p\pi_n^*H^0({\cal
{A}}_g^P(n) ,\o^{\otimes k})$. Let $d=deg(\pi^{}_{T,B})$. Then
$df\pm g$ are eigenforms of $\pi^{}_{T,B_*} v^{\sigma
*}_p\pi^{*}_{T,B}$ with eigenvalue $\pm d$. Now, the truth of
proposition is a consequence of the first stage and the following
\begin{lem}
Any nonzero Siegel modular form on ${\cal {A}}_g^{B,P}(p,n)$ which
is an eigenform of $\pi^{}_{T,B_*} v^{\sigma *}_p\pi^{*}_{T,B}$
for a non-zero $\sigma$ in $W_G$ is pull back of a Siegel form on
${\cal {A}}_g^{P',P}(p,n)$ where $P'$ is the parabolic subgroup
generated by $B$ and $\sigma B\sigma^{-1}$.
\end{lem}
\textbf {Proof.} This is a consequence of propositions 1.1. and
1.2. $\square$

In the final stage we show that the above-mentioned $g!2^g$
copies of level-$n$ Siegel forms inside the space of Siegel forms
of level $pn$ are indeed $p$-old forms of level-$pn$ and they are
linearly independent.
\begin{prop}
The Atkin-Lehner correspondences $C^{\sigma}_B(p)\pi^*_n$ commute
with all Hecke correspondences which generate the local Hecke
algebras $H_q$ for $q$ relatively prime to $pn$.
\end{prop}
\textbf {Proof.} The action of Atkin-Lehner correspondences can
be interpreted in terms of the $p$-torsion of abelian varieties
representing points of the moduli-space. By geometric base-change
one can see that such an action commutes with those interpreted
in terms of the $q$-torsion points.$\square$ \\\textbf{Proof of
the main theorem.} The $g!2^g$ copies of level-$n$ Siegel forms
induced by Atkin-Lehner correspondences are contained in the
space of $p$-old forms by previous proposition. Recall that in
this paper Siegel modular forms with eigenvalues repeated from
lover levels are considered to be old. The space of Siegel forms
on ${\cal {A}}_g^{B,P}(p,n)$ is finite-dimensional and has a basis
of all prime to $pn$ Hecke eigenforms. So is the case for any of
the $g!2^g$ Atkin-Lehner copies, by previous proposition. Fix a
basis for the space of Siegel forms of level $pn$ whose elements
are eigenforms of all prime to $pn$ local Hecke algebras. Suppose
the $g!2^g$ Atkin-Lehner copies are linearly dependent. It means
that a non-zero eigenform $f$ of all prime to $pn$ local Hecke
algebras is generated by some basis elements of the Atkin-Lehner
copies which have the same eigenvalues. Consider the vector space
$V$ of all Siegel forms with the same Hecke eigenvalues $f$ and
consider the $g!2^g$ Atkin-Lehner copies in it. $V$ is invariant
under the action of all the $g!2^g$ correspondences
$D^{\sigma}_B(p)=\pi^{}_{T,B_*} v^{\sigma *}_p\pi^{*}_{T,B}$ for
$\sigma\in W_G$, because these correspondences commute with the
Atkin-Lehner copies. By this symmetry, $df+\pi^{}_{T,P_*}
v^{\sigma *}_p\pi^{*}_{T,P} f$ is also a Hecke eigenform in $V$
with the same Hecke eigenvalues as $f$. There exists a $\sigma\in
W_G$ which gives a nonzero Siegel form $df+\pi^{}_{T,B_*}
v^{\sigma *}_p\pi^{*}_{T,B} f$ generated by the Atkin-Lehner
copies. But such a vector is a Siegel form which can be pulled
back from a lower Siegel space by proposition 1.2.. This
contradicts linear dependency.$\square$

Having constructed the Atkin-Lehner copies of $p$-old forms on
${\cal {A}}_g^{B,P}(p,n)$ one can get the $2^g$ copies of $p$-old
forms on ${\cal {A}}_g^P(pn)$ by pushing forward all the $g!2^g$
$p$-old copies down to ${\cal {A}}_g^P(pn)$.

\begin{thm}
The linear subspaces of $H^0({\cal {A}}_g^P(pn) ,\o^{\otimes k})$
generated by Atkin-Lehner correspondences $C^{\sigma}_P(p)\pi^*_n$
for $\sigma$ varying in $W_G$ where $C^{\sigma}_P(p)$ is defined
by
$$
\sum \pi^{}_{T,P*} v^{\eta\sigma *}_p\pi^{*}_{T,B}w^{B*}_p
$$
where the sum ranges over $\eta\in S_g\subset W_G$ give $2^g$
linearly independent copies of $H^0({\cal {A}}_g^P(n)
,\o^{\otimes k})$ inside $p$-old forms of level $pn$ living on
${\cal {A}}_g^P(pn)$.
\end{thm}
\textbf {Proof.} This is a direct consequence of the main theorem
and proposition 1.5. $\square$

\section{Algebraic equations for towers of Siegel spaces}

\subsection{Construction of ${\cal {A}}_g^{P}(p^k)$ from ${\cal {A}}_g^{P}(p)$'s and ${\cal {A}}_g^{P}(p^2)$'s}

In this section, we follow Elkies, who did the genus one case [Elk].
Fix a prime $p>1$.  For positive $k$,
the Siegel moduli space ${\cal {A}}_g^{P}(p^k)$
parametrizes principally polarized abelian varieties 
with a cyclic $p^k$-isogeny, or equivalently sequences of
$p$-isogenies
$$
A_0 \rightarrow A_1 \rightarrow A_2 \rightarrow \cdots \rightarrow A_k
$$
such that the composite isogeny $A_{j-1} \rightarrow A_{j+1}$ of degree $p^{2g}$
is cyclic for each $j$ with $0<j<k$.  Thus for each $m=0,1,\ldots,k$
there are $k+1-m$ maps $\pi_j:{\cal {A}}_g^{P}(p^k) \rightarrow {\cal {A}}_g^{P}(p^m)$ obtained by
extracting for some $j=0,1,\ldots,k-m$ the cyclic $p^m$-isogeny
$A_j \rightarrow A_{j+m}$ from the above sequence.  
In particular we have a tower of maps
$$
{\cal {A}}_g^{P}(p^k) \rightarrow {\cal {A}}_g^{P}(p^{k-1}) \rightarrow {\cal {A}}_g^{P}(p^{k-2}) \rightarrow
 \cdots \rightarrow {\cal {A}}_g^{P}(p^2) \rightarrow {\cal {A}}_g^{P}(p)
$$
each map being of degree $p^g$.  Each ${\cal {A}}_g^{P}(p^k)$ also has an
Atkin-Lehner involution $w_p=w_p^{(k)}$, taking a cyclic
$p^k$-isogeny to its dual isogeny, and the above sequence 
to the sequence
$$
A_k \rightarrow \cdots \rightarrow A_2 \rightarrow A_1 \rightarrow A_0
$$
of dual isogenies.  We thus have
$$
w_p^{(m)} \circ \pi_j = \pi_{k-m-j} \circ w_p^{(k)},
$$
where $\pi_j, \pi_{k-m-j}$ are our $j$th and $(k-m-j)$th
maps from ${\cal {A}}_g^{P}(p^k)$ to ${\cal {A}}_g^{P}(p^m)$.

Now the explicit formulas for
${\cal {A}}_g^{P}(p),{\cal {A}}_g^{P}(p^2)$, together with the involutions
$w_p^{(1)},w_p^{(2)}$ of these moduli spaces and the map
$\pi_0:{\cal {A}}_g^{P}(p^2)\rightarrow {\cal {A}}_g^{P}(p)$ between them, suffice to exhibit
the entire tower explicitly:

begin{prop}  
For $k\geq 2$ the product map
$$
\pi = \pi_0 \times \pi_1 \times \pi_2 \times \cdots \times \pi_{k-2}:
{\cal {A}}_g^{P}(p^k) \rightarrow ({\cal {A}}_g^{P}(p^2))^{k-1}
$$
is a 1:1 map from ${\cal {A}}_g^{P}(p^k)$ to the set of
$(P_1,P_2,\ldots,P_{k-1}) \in \bigl({\cal {A}}_g^{P}(p^2)\bigr)^{k-1}$ such that
$$
\pi_0 \bigl( w_p^{(2)} (P_j) \bigr) =
w_p^{(1)} \bigl( \pi_0(P_{j+1}) \bigr)
$$
for each $j=1,2,\ldots,k-2$.
end{prop}

Informally speaking, we get from ${\cal {A}}_g^{P}(p^2)$ up to ${\cal {A}}_g^{P}(p^k)$ by
iterating $k-2$ times the involution $w_p^{(2)}$ composed with the
''$p$-valued involution'' $\pi_0^{-1} w_p^{(1)} \pi_0$.
Of course the maps $\pi_j: {\cal {A}}_g^{P}(p^k) \rightarrow {\cal {A}}_g^{P}(p^m)$ (for $m\geq 2$)
are then simply
$$
(P_1,\ldots,P_{k-1}) \mapsto (P_{j+1},\ldots,P_{j+m-1}),
$$
and the involution $w_p^{(k)}$ is
$$
(P_1,P_2,\ldots,P_{k-2},P_{k-1}) \leftrightarrow
(w_p^{(2)}P_{k-1},w_p^{(2)}P_{k-2},\ldots,w_p^{(2)}P_2,w_p^{(2)}P_1),
$$
i.e.\ reversing the order of $P_1,\ldots,P_{k-1}$ and applying
$w_p^{(2)}$ to each coordinate.

\textbf {Proof.} 
It is clear that the map is 1:1 to its image, because
a sequence of $p$-isogenies is determined by the
$p^2$-isogenies $A_{j-1} \rightarrow A_{j+1}$ parametrized by the
$j$th coordinate of~$\pi$ ($0<j<k$).  Now $(P_1,\ldots,P_{k-1})$
is in the image of~$\pi$ if and only if the $p^2$-isogenies
parametrized by $P_1,\ldots,P_{k-1}$, regarded as sequences
$A_0^j \rightarrow A_1^j \rightarrow A_2^j$ of $p$-isogenies, fit together
to form a sequence with $A_i^j=A_{i+j}$,
i.e.\ if and only if the isogenies $A_1^j \rightarrow A_2^j$ and
$A_0^{j+1} \rightarrow A_1^{j+1}$ coincide for each $j=1,2,\ldots,k-2$.
But these isogenies are represented by the points $\pi_1(P_j)$
and $\pi_0(P_{j+1})$ on ${\cal {A}}_g^{P}(p)$.
Thus the necessary and sufficient condition is that
$$
\pi_1(P_j) = \pi_0(P_{j+1})
$$
for each $j=1,2,\ldots,k-2$; applying $w_p^{(1)}$ to both sides,
and then commutativity of the Atkin-Lehner involutions to $w_p^{(1)} \bigl( \pi_1(P_j) \bigr)$,
then yields the equivalent form of what we seek. $\square$

\subsection{Construction of ${\cal {A}}_g^{P}(n)$ from ${\cal {A}}_g^{P}(p^k)$'s}

\begin{prop}
 Let $n=p_1^{\alpha_1} \dots p_k^{\alpha_k}$ be the prime decomposition of  an integer $n\geq 2$
then if $\pi_j: {\cal {A}}_g^{P}(p_j^{\alpha_j}) \rightarrow {\cal {A}}_g(1)$ denotes the natural 
projection forgetting the level structure for $j=1,...,k$, then ${\cal {A}}_g^{P}(n)$ is nothing but 
the fiber product of ${\cal {A}}_g^{P}(p_j^{\alpha_j})$'s over ${\cal {A}}_g(1)$ via $\pi_j$.
\end{prop}

In fact, the following is true:
\begin{prop}
 Let $n$ and $m$ relatively prime natural numbers.
Let $\pi': {\cal {A}}_g^{P}(n) \rightarrow {\cal {A}}_g(1)$ and 
$\pi'': {\cal {A}}_g^{P}(m) \rightarrow {\cal {A}}_g(1)$
denote the natural 
projections forgetting the level structures, then 
$$
{\cal {A}}_g^{P}(nm)={\cal {A}}_g^{P}(n)\times_{{\cal {A}}_g(1)} {\cal {A}}_g^{P}(m)
$$ 
where $\times$ is the fiber product over ${\cal {A}}_g(1)$ via $\pi'$ and $\pi''$.
\end{prop}
\textbf {Proof.} 
Considering the moduli interpretation of the above mentioned moduli spaces, the fact that 
$n$-level atructures are independent of $m$-level structures for $(n,m)=1$ and fixing these two we can get an
$mn$-level structure gives a one-to-one map from the right hand side to the left.
$\square$

By the previous section, once we have algebraic equations for ${\cal {A}}_g^{P}(p)$ and ${\cal {A}}_g^{P}(p^2)$
and morphisms over ${\cal {A}}_g(1)$ we get algebraic equations for ${\cal {A}}_g^{P}(p^k)$ over ${\cal {A}}_g(1)$
for all primes $p$ and this will suffice to get algebraic equations for all ${\cal {A}}_g^{P}(n)$ as we desire.

\subsection{Compactification of Siegel moduli spaces}

The space of Siegel modular forms can also be formulated in the
language of schemes. Let $S$ be a base scheme. A modular form $f$
of weight $k$ is a rule which assigns to each principally
polarized abelian variety $(A/S,\l)$ a section $f(A/S,\l)$ of
$\o_{A/S}^{\otimes k}$ over $S$ depending only on the isomorphism
class of $(A/S,\l)$ commuting with arbitrary base change. Here
$\o_{A/S}$ is the top wedge of tangent bundle at origin of $A$
over $S$.

To define Siegel modular forms of higher level, one should equip
principally polarized abelian varieties with level structures. Let
$\z_n$ denote an $n$-th root of unity where $n\geq 3$. On a
principally polarized abelian scheme $(A,\l)$ over $\Sn$ of
relative dimension $g$ we define a symplectic principal level-$n$
structure to be a symplectic isomorphism
$\a:A[n]\to(\Z/n\Z)^{2g}$ where $(\Z/n\Z)^{2g}$ is equipped with
the standard non-degenerate skew-symmetric pairing
$$
<,>:(\Z/n\Z)^{2g} \times (\Z/n\Z)^{2g} \to \Z/n\Z
$$
$$
<(u,v) , (z,w)> \mapsto u.w^t-v.z^t
$$

Let $S$ be a scheme over $\Sn$. The moduli scheme classifying the
principally polarized abelian schemes over $S$ together with a
symplectic principal level-$n$ structure is a scheme over $S$
and  will be denoted by $\Agn$. The moduli scheme $\Agn$ over $S$
can be constructed from $\Agn$ over $\Sn$ by base change.

$Sp(2g,\Z/n\Z)$ acts as a group of symmetries on $\Agn$ by acting
on level structures. We will recognize these moduli spaces and
their equivariant quotients under the action of subgroups of
$Sp(2g,\Z/n\Z)$ as Siegel spaces. We restrict our attention to
Siegel spaces over $\Sn$. $\Agn$ is connected and smooth over
$\Sn$. The condition $n \geq 3$ is to guarantee that we get a
moduli scheme, instead of getting only a moduli stack. The
natural morphism $\Agn \to \Agm$ where $m$,$n$ are positive
integers $\geq 3$ and $m|n$ is a finite and etale morphism over
$\Sn$.

Let $\Bgn$ denote the universal abelian variety over $\Agn$. The
Hodge bundle $\o$ is defined to be the pull back via the zero
section $i_0:\Agn \to \Bgn$ of the line bundle
$\wedge^{top}\Omega_{\Bgn/\Agn}$. The Hodge bundle is an ample
invertible sheaf on $\Agn$ and can be naturally extended to a
bundle $\o$ on $\Acgn$. We could define the minimal
compactification $\Acgn$ by the formula
$$
\Acgn =proj(\oplus_{k \geq 0}H^0(\Agn ,\o^{\otimes k})).
$$
The graded ring above is regarded as a $\Z[\z_n,1/n]$-algebra. The
scheme $\Acgn$ is equipped with a stratification by locally
closed subschemes which are geometrically normal and flat over
$\Sn$. Each of these strata is canonically isomorphic to a moduli
space $A_{i}(n)$ for some $i$ between $0$ and $g$. The map $\Agn
\to \Agm$ can be extended uniquely to $\Acgn \to \Acgm$ for
$m|n$. These maps when restricted to strata, induce the
corresponding natural maps between lower genera Siegel spaces
$A_i(n) \to A_i(m)$. The action $Sp(2g,\Z/n\Z)$ on $\Agn$
naturally extends to an action on the compactified Siegel space
$\Acgn$. This action is compatible with the maps $\Acgn \to
\Acgm$ for $m|n$.

Let $K(n)$ denote the subgroup of $Sp(2g,\Z/n\Z)$ fixing the $g$
first $(\Z/n\Z)$-basis elements of $(\Z/n\Z)^{\oplus 2g}$ on
which $Sp(2g,\Z/n\Z)$ acts. Since $\Acgn$ is a projective scheme,
we can define the quotient projective schemes $\Aocgn$ to be the
geometric quotient of $\Acgn$ by $K(n)$. This quotient provides
us with a compactification of $\Aogn$ which is the moduli scheme
of principally polarized abelian schemes $(A,\l)$ over $\Sn$,
together with $g$ elements in $A[n]$ generating a symplectic
subgroup isomorphic to $(\Z/n\Z)^g$. Again we have natural maps
$\Aocgn \to \Aocgm$ for $m|n$.

We define the Hodge bundle $\o$ on $\Aocgn$ to be the quotient of
the Hodge bundle $\o$ on $\Acgn$ under the action of the
corresponding subgroup $K(n)$ of $Sp(2g,\Z/n\Z)$. This is
possible because the line bundle $\o$ on the space $\Acgn$ is
$Sp(2g,\Z/n\Z)$-linearizable. A Siegel modular form of weight $k$
and full level $n$ is a global section of $\o^k$ on $\Acgn$. Over
the complex numbers, this corresponds to a Siegel modular form of
weight $k$ with respect to $\Gamma(n)$. In this paper, by a
Siegel modular form of weight $k$ and level $n$ we mean a global
section of $\o^k$ on $\Aocgn$. This corresponds to the congruence
subgroup $\Gamma_0(n)$.

\subsection{Algebraic equations for ${\cal {A}}_2^*(2)$}

In this section, we follow Lee and Weintraub [Lee-Wei1-5]. 
For construction of compactifications of ${\cal {A}}_2(1)$ look at 
[Lee-Wei3]. The compactification we have built in the previous section
is called the Satake compactification which is a projective algebraic variety with severe 
singularities. It would be also handy to introduce algebraic equations for smooth compactifications
as was constructed by Igusa for $g=2$ and generalized by Mumford and his collaborators 
to general genus and extended to schemes by Falting and Chai [Fal-Chai] 
which is called the toroidal compactification. 

The Siegel moduli space ${\cal {A}}_2^*(2)$ is related another variety which appeared
in the work of Deligne and Mostow [Del-Mos], constructed by means of Mumford’s geometric invariant theory. 

Let $S$ denote the set $\{ 1,2,3,4,5,6\}$, and let $\Bbb P^5$ denote the space of functions of $S$ 
to $\Bbb P^1 = \Bbb P^1(\Bbb C)$. There is a natural action of $PGL_2(\Bbb C)$ on the space $\Bbb P^1$ induced by the linear 
fractional transformation of $PGL_2(\Bbb C)$ on $\Bbb P^1$. The subspace of injective functions can be 
identified with $(\Bbb P^1)^6-\Delta$, and its quotient with the moduli space $\cal M$ of nonsingular curves 
with level 2 structure. By a stable point (resp. semi-stable point) of $\Bbb P^5$, we mean a point 
with the property that no more than two (resp. three) elements in $S$ have the same 
image. The group $PGL_2$ operates freely on the subspace of stable points, and its 
quotient space ${\cal {Q}}_{st}$, is a quasi-projective variety. To compactify ${\cal {Q}}_{st}$, we consider the space of 
semi-stable points. Define an equivalence relation for which two stable points are equivalent if and only if 
they have the same $PGL_2$-orbit and if two points are semi-stable but not stable 
they are equvalent if they induce the same partition of $S$ into two sets of three elements $S_1$ and $S_2$ such
that each functaion separates them and on constant on one of them. 
The quotient space ${\cal {Q}}_{sst}$ of semi-stable points module this relation is a projective variety, and 
contains ${\cal {Q}}_{st}$ as a Zariski open set. In fact, ${\cal {Q}}_{cusp}={\cal {Q}}_{sst}-{\cal {Q}}_{st}$, consists of ten isolated singular 
points. To desingularize this variety, we blow up these points and obtain a nonsingular variety ${\cal {Q}}^{\sim}_{sst}$
\begin{prop}
${\cal {Q}}^{\sim}_{sst}$ is isomorphic to the Igusa 
compactification ${\cal {A}}_2^{\sim}(2)$
\end{prop}

The Igusa compactification ${\cal {A}}_2^{\sim}(2)$ of ${\cal {A}}_2(2)$ may be constructed by desingularizing, or blowing 
up, the Satake compactificationx of${\cal {A}}_2^{*}(2)$. 
Lee and Weintraub construct a birational transformation $f:{\cal {A}}_2^{*}(2) \rightarrow {\cal {Q}}_{sst}$.
To begin, they identify ${\cal {Q}}_{sst}$, with a classical object, Segre’s cubic threefold. 
From this it follows that ${\cal {Q}}_{sst}$ is isomorphic to the threefold in 
$\Bbb P^5$ defined by the homogeneous equations 
$$
\sum_{i=1}^6 x_i=0
$$
$$
\sum_{i=1}^6 x^3_i=0
$$
known as Segre’s cubic threefold.
Since Segre’s time it has been known that considering the dual hypersurface to the symmetric quartic threefold defined by 
the equations 
$$
\sum_{i=1}^6 x_i=0
$$
$$
\left(\sum_{i=1}^6 x^4_i\right)-\left(\sum_{i=1}^6 x^2_i\right)^2=0
$$
yields Segre’s cubic threefold. van der Geer [vdG] has shown that the quartic threefold defined above can be identified 
with the Satake compactification ${\cal {A}}_2^{*}(2)$, so one obtains a birational transformation 
$$
\tilde {f} :{\cal {A}}_2^{*}(2) \dashrightarrow {\cal {Q}}_{sst}
$$
Alternately, we may consider the Igusa compactification ${\cal {A}}_2^{\sim}(2)$. Note that the Satake 
compactification ${\cal {A}}_2^{*}(2)$ is a hypersurface in $\Bbb P^5$ defined by a single function $F(x_1, . . . , x_5) = 0.$
The derivatives $\partial F/\partial x_i$, on the one hand define the coordinate functions of the projective dual, 
and on the other hand generate the ideal $I = ( \partial F /\partial x_1, . . . , \partial F/\partial x_5 )$ that defines the 
boundary components $\partial$. Since ${\cal {A}}_2^{\sim}(2)$ is defined by blowing up ${\cal {A}}_2^{*}(2)$ along $\partial$, 
it follows that $\tilde {f}$ lifts to a morphism 
from ${\cal {A}}_2^{\sim}(1)$ to the projective dual of ${\cal {Q}}_{sst}$. one can blow down the ten components of the Humbert surface in
${\cal {A}}_2^{\sim}(1)$ to 
points to get a complex analytic space ${\cal {A}}_2^{*}(2)$. Then, from the definitions, we have the mapping $\hat {f}$ in 
the following diagram: 
$$
\begin{array}{ccccc}
\hat {f}:{\cal {A}}_2^{\sim}(2) & 
\longrightarrow & {\cal {Q}}^{\sim}_{sst}  \\
\downarrow & & \downarrow \\
\tilde {f}:{\cal {A}}_2^{*}(2) & 
\longrightarrow & {\cal {Q}}^{*}_{sst}  
\end{array}
$$

It can be shown that $\hat {f}$ and $f^{\sim}$ are isomorphisms.

By division by the two-level structure on ${\cal {Q}}_{sst}$, we can find the ring of invariants of ${\cal {A}}_2^{*}(1)$
and the relative morphism of ${\cal {A}}_2^{*}(2)$ over ${\cal {A}}_2^{*}(1)$.

\subsection{Algebraic equations for ${\cal {A}}_2^{*}(3)$}

In this section, we follow Hoffman and Weintraub [Hof-Wei].
Here, one uses  a variety $\cal B$ defined over $\Bbb Q(\sqrt{-3})$
such that $\cal B$ over $\Bbb C$ is isomorphic with ${\cal {A}}_2^{*}(3)$. 
Felix Klein initiated the study of the moduli spaces of genus 2 curves and the coverings defined by "Stufe".
Two of his students, H. Burkhardt and H. Maschke, took up the case where
Stufe = 3. Burkhardt managed to write down an explicit equation for this moduli
space. The general idea is this: Consider the 9 thetanullwerte
$$
X_{\alpha ,\beta}=\theta \left(   
\begin{array}{ccc} 0 & 0 \\ \alpha & \beta \end{array}
\right) (\tau ,0),
~~~~~~\alpha \in  1/3\Bbb Z/\Bbb Z  ~~\beta \in 1/3\Bbb Z/\Bbb Z
$$
These 9 values have the property that as $\tau \mapsto \gamma .\tau$ with $\gamma \in P\Gamma_2(1)$ they undergo a
linear transformation, which is the identity up to scalar multiples for $\gamma \in P\Gamma_2(3)$.
In other words, we have a projective representation of the finite simple group of
order 25920, $G = P\Gamma_2(1) / P\Gamma_2(3)= PSp(4;\Bbb F_3)$. This representation splits into
two invariant subspaces of dimension 4 and 5 respectively, the spaces of
$$
Z_{\alpha ,\beta}=(X_{\alpha ,\beta}-X_{-\alpha ,-\beta})/2 ~~~
Y_{\alpha ,\beta}=(X_{\alpha ,\beta}+X_{-\alpha ,-\beta})/2
$$
Maschke studied the action of $G$ on the $Z$'s, Burkhardt studied the action on the
$Y$'s, and both managed to find the ring of $G$-invariant forms in their respective
cases. Let
$$
Y_0= Y_{0,0}, ~~ 2Y_1= Y_{1/3,0}, ~~2Y_2= Y_{0,1/3},~~2Y_3= Y_{1/3,1/3}, ~~2Y_4= Y_{1/3,2/3}.
$$
Burkhardt found the invariant form of degree 4:
$$
J_4=Y_0^4-Y_0(Y_1^3+Y_2^3+Y_3^3+Y_4^3)+3Y_1Y_2Y_3Y_4.
$$
\begin{prop} 
 Let ${\cal {B}}_0 \subset \Bbb P^4$ be the quartic hypersurface defined by $J_4=0$
There is an isomorphism between a Zariski open subset of ${\cal {A}}_2(3)$
and a Zariski open subset of ${\cal {B}}_0$.
Let ${\cal {B}}$ be the variety obtained by resolving the 45 nodes on ${\cal {B}}_0$. The the map
above extends to an isomorphism with the Igusa compactification: ${\cal {B}}_0\simeq {\cal {A}}_2^{*}(3)$.
\end{prop}

In this form the proposition was first proved by van der Geer [vdG2], who asserted some
thing stronger, namely that these results were true for the corresponding schemes
over $\Bbb Z[1/3; \epsilon]$, where $\epsilon$ is a primitive cube root of unity (the existence of a model of
${\cal {A}}_2^{*}(3)$ over that ring being a consequence of Faltings' theory [Cha-Fal]). 

Since we have the group action explicitly, we can find the ring of invariants of ${\cal {A}}_2^{*}(1)$
and the relative morphism of ${\cal {A}}_2^{P*}(3)$ over ${\cal {A}}_2^{*}(1)$.

\subsection{Algebraic equations for ${\cal {A}}_2^{*}(4)$}

In this section, we follow van Geeman and Nygaard [vGe-Nyg] and its exposition by  Okazaki and Yamauchi [Oka-Yam].
Let $\mathcal{A}_2(2,4,8)$ be the moduli space of abelian surfaces 
with some level structure which has been studied by van Geemen and Nygaard. 
It is the quotient space of the Siegel upper half plane  of degree 2  
by the arithmetic subgroup $\Gamma(2,4,8)$ of the symplectic group $Sp_4(\Z)$. 
This congruence subgroup $\Gamma(2,4,8)$ is contained in the principal congruence subgroup 
$\Gamma(4):=\{\gamma\in Sp_4(\Z)\ |\ \gamma\equiv 1_4\ {\rm mod}\ 4 \}$. 
$\mathcal{A}(2,4,8)$ is a quasi-projective smooth threefold. 
By \cite{vG-N} we have the projective model ${\mathcal{A}}^{*}_2(2,4,8)$ the Satake compactification of $\mathcal{A}_2(2,4,8)$ which 
is defined over $\Q$ in $\Bbb P^{13}$ as follows:
$$
\begin{array}{rl}
Y^2_0&=Q_0(X_0,X_1,X_2,X_3):=X^2_0+X^2_1+X^2_2+X^2_3\\
Y^2_1&=Q_1(X_0,X_1,X_2,X_3):=X^2_0-X^2_1+X^2_2-X^2_3\\
Y^2_2&=Q_2(X_0,X_1,X_2,X_3):=X^2_0+X^2_1-X^2_2-X^2_3\\
Y^2_3&=Q_3(X_0,X_1,X_2,X_3):=X^2_0-X^2_1-X^2_2+X^2_3\\
Y^2_4&=Q_4(X_0,X_1,X_2,X_3):=2(X_0X_1+X_2X_3)\\
Y^2_5&=Q_5(X_0,X_1,X_2,X_3):=2(X_0X_2+X_1X_3)\\
Y^2_6&=Q_6(X_0,X_1,X_2,X_3):=2(X_0X_3+X_1X_2)\\
Y^2_7&=Q_7(X_0,X_1,X_2,X_3):=2(X_0X_1-X_2X_3)\\
Y^2_8&=Q_8(X_0,X_1,X_2,X_3):=2(X_0X_2-X_1X_3)\\
Y^2_9&=Q_9(X_0,X_1,X_2,X_3):=2(X_0X_3-X_1X_2).
\end{array}
$$
Since we have the group action of $Sp_4(\Z)$ explicitly, we can find the ring of invariants of ${\cal {A}}_2^{*}(4)$ and ${\cal {A}}_2^{*}(1)$
and the relative morphism of ${\cal {A}}_2^{P*}(4)$ over ${\cal {A}}_2^{*}(1)$.

\subsection{Algebraic equations for ${\cal {A}}_2^{P*}(2^k)$ and ${\cal {A}}_2^{P*}(3.2^k)$}

By section 2.1. we can get algebraic equations for ${\cal {A}}_2^{P*}(2^k)$ using
algebraic equations for ${\cal {A}}_2^{P*}(4)$ and 
${\cal {A}}_2^{P*}(2)$ over ${\cal {A}}_2^{P*}(1)$. By section 2.2. using algebraic equations for 
${\cal {A}}_2^{P*}(2^k)$ and ${\cal {A}}_2^{P*}(3)$ over ${\cal {A}}_2^{P*}(1)$ one can get algebraic equations for 
${\cal {A}}_2^{P*}(3.2^k)$.

\subsection*{Acknowledgements}

For this research we have benefited from conversations with
J. de Jong, A. Genestier,
G. Pappas, A. Rajaei, P. Sarnak, C. Skinner, R. Takloo-bighash, R. Taylor, J. Tilouine and A.
Wiles. We shall thank A. Genestier
and J. Tilouine heartily for many useful comments and questions on an early version
which lead us to correction of many mistakes.


\thispagestyle{plain}        
\addcontentsline{toc}{chapter}{Bibliography}
\bibliographystyle{amsplain}

Sharif University of Technology, e-mail: rastegar@sharif.ir
\\ Princeton University, e-mail:
rastegar@math.princeton.edu

\end{document}